\documentclass[12pt]{article}
\usepackage{amsmath}
\usepackage{amsfonts}
\usepackage{amssymb}
\usepackage{mathtools}
\usepackage{amscd}
\newtheorem{theorem}{Theorem}[section]
\newtheorem{remark}[theorem]{Remark}

\newtheorem{example}[theorem]{Example}
\newtheorem{corollary}[theorem]{Corollary}
\newtheorem{proposition}[theorem]{Proposition}

\begin{document}

\begin{center}
{\Large \bf On direct and inverse limits of Minkowski's balls, domains, and their critical lattices}
\end{center}

\begin{center}
{\bf N.M. Glazunov} 
 \end{center}

\begin{center}
{\rm Glushkov Institute of Cybernetics NASU, Kiev, } \\
{\rm Institute of Mathematics and Informatics Bulgarian Academy of Sciences }\\
{\rm  Email:} {\it glanm@yahoo.com }
\end{center} 

\bigskip

\footnote{ The author was supported   Simons grant 992227.}

\bigskip
${\mathcal  \; 2020\;  Mathematics \; Subject  \; Classification}$:  11H06, 52C05

\bigskip


${\mathcal  \;  Key\; words \; and \; phrases}$:  Lattice packing,  Minkowski's ball, Minkowski's domain, critical lattice, optimal lattice packing, direct system,   direct limit. inverse system, inverse limit
\bigskip

{\bf Abstract.} We construct direct and inverse systems of Minkowski's  balls and domains, direct and inverse systems of their critical lattices and calculate their direct and inverse limits.\\

\section{Introduction}
In this paper we continue the investigations begun in (\cite{Gl5}) 
of optimal packing of  Minkowski's balls in the plane ${\mathbb R}^2$. 
We construct direct and inverse Minkowski systems (including the cases of Minkowski, Davis and Chebyshev-Cohn) of balls and domains, direct and inverse systems of their critical lattices, and calculate their direct and inverse limits.
In section 2 we present  expressions for critical determinants 
$\Delta(D_p)$ and corresponding critical  lattices:
Section 3 deals with Minkowski, Davis and Chebyshev-Cohn balls and domains.
In section 4 we construct direct systems and direct limits and in section 5 
we construct inverse systems and inverse limits.
For simplicity we present results for Minkowski's (or unit $l_p$ norms) balls
 and Minkowski's domains and their critical lattices and we leave the reader 
the  opportunity to specify (on the basis of Theorem (\ref{mt}) and the 
definitions of Section 3) these results for balls, domains, and lattices of 
Minkowski, Davis and Chebyshev-Cohn.
Direct systems and direct limits (as well as inverse systems and inverse limits) were defined by Pontryagin \cite{Pontryagin}.
\v{C}ech has used  these notions  implicitly under construction of his cohomology theory. 
Much progress has been made in this theory after the publishing of the book \cite{stei} through  the work by. J.P. Serre,
A. Grothendieck, P. Deligne, J. Tate and their followers. 
To the best of our knowledge the direct and inverse systems and their  limits for optimal packing problems are introduced for the first time in the present paper.

\section{Minkowski's balls and their critical lattices}

By general ( Minkowski's) balls we  mean (two-dimensional) balls (unit $l_p$ norms) in ${\mathbb R}^2$ of the form
  \begin{equation}
   \label{eq1}
   D_p: \;  |x|^p + |y|^p \le 1, \; p \ge 1. 
  \end{equation}
 \begin{remark}
  Minkowski's balls  are semialgebraic sets. A semialgebraic set is a finite union of sets defined by polynomial equalities and polynomial inequalities. A semialgebraic function is a function with a semialgebraic graph. 
For $p = 2d, d-$natural, the boundary of   Minkowski's balls are algebraic curves.
 \end{remark}
 
 Denote by $V(D_p)$ the volume (area) of $D_p$.

\begin{proposition} [\cite{Mi:DA}]
 The volume of   Minkowski's ball $D_p$ is equal $4 \frac{(\Gamma(1 + \frac{1}{p}))^2}{\Gamma(1 + \frac{2}{p})}$.
\end{proposition}

 

From the proof of Minkowski's conjecture  \cite{Mi:DA,M:LP,D:NC,Co:MC,W:MC,GGM:PM} in notations \cite{GGM:PM,Gl5}
   we have next expressions for critical determinants $\Delta(D_p)$ and corresponding critical  lattices:  
   \begin{theorem}
   \label{mt}
\begin{enumerate} 

\item  $\Delta(D_p) = {\Delta^{(0)}_p} = \Delta(p, {\sigma_p}) =  \frac{1}{2}{\sigma}_{p}, \;  2 \le p \le p_{0};  $

\item $ {\sigma}_{p} = (2^p - 1)^{1/p},$

\item  $\Delta(D_p) ={\Delta^{(1)}_p}  = \Delta(p,1) = 4^{-\frac{1}{p}}\frac{1 +\tau_p }{1 - \tau_p}, \; 1 \le p \le 2, \; p \ge p_{0}, $

\item  $2(1 - \tau_p)^p = 1 + \tau_p^p,  \;  0 \le \tau_p < 1,$  
\end{enumerate}
here $p_{0}$ is a real number that is defined unique by conditions
$\Delta(p_{0},\sigma_p) = \Delta(p_{0},1),  \;
2,57 < p_{0}  < 2,58, \; p_0  \approx 2.5725 $\\
  For their critical lattices respectively  $\Lambda_{p}^{(0)},\; \Lambda_{p}^{(1)}$ next conditions satisfy:   $\Lambda_{p}^{(0)}$ and 
 $\Lambda_{p}^{(1)}$  are  two $D_p$-admissible lattices each of which contains
three pairs of points on the boundary of $D_p$  with the
property that 
 $(1,0) \in \Lambda_{p}^{(0)},$ 
$(-2^{-1/p},2^{-1/p}) \in \Lambda_{p}^{(1)},$
\end{theorem}



\begin{example}
 Lattices $\Lambda_{p}^{(0)}$ are two-dimensional lattices in ${\mathbb R}^2$ spanned by the vectors
 \begin{itemize}
 \item[]  $\lambda^{(1)} = (1, 0),$
   \item[] $\lambda^{(2)} = (\frac{1}{2}, \frac{1}{2} \sigma_p).$
   \end{itemize}
  
    The lattice $\Lambda_{2}^{(1)}$ is a two-dimensional lattice in ${\mathbb R}^2$ spanned by the vectors
 \begin{itemize}
 \item[]  $\lambda^{(1)} = (-2^{-1/2},2^{-1/2}),$
   \item[] $\lambda^{(2)} = (\frac{\sqrt 6 - \sqrt 2}{4}, \frac{\sqrt 6 + \sqrt 2}{4}).$
   \end{itemize}
   \end{example}

\section{ Balls and Domains}

We consider balls of the form
\[
   D_p: \;  |x|^p + |y|^p \le 1, \; p \ge 1,
\]
and call such balls with $1<p<2$ {\it Minkowski balls}.
Continuing this, we consider the following classes of balls and circles.
\begin{itemize}
  \item {\it Davis balls}: $|x|^p + |y|^p \le 1$  for $p_{0} > p \ge 2$;
  \item {\it Chebyshev-Cohn balls}: $|x|^p + |y|^p \le 1$  for $ p \ge p_{0}$;
  \end{itemize}
  
    Let $D$ be a fixed bounded symmetric about origin  convex body ({\it centrally symmetric convex body} for short) with volume
$V(D)$.
\begin{proposition} \cite{Cassels}.
\label{p1}
 If $D$ is symmetric about the origin and convex, then $2D$ is convex and symmetric
 about the origin.
 \end{proposition}  
 \begin{corollary}
 \label{cor1}
 Let $m$ be integer $m \ge 0$ and $n$ be natural greater $m$. 
If $2^m D$ centrally symmetric convex body then $2^n D$ is again centrally symmetrc convex body.
\end{corollary}
{\bf Proof.} Induction. \\

  We consider the following classes of balls (see above) and domains.
\begin{itemize}
  \item  {\it  Minkowski domains}:    $2^m D_p$,  integer $m \ge 1$, for $1 \le p<2$;
\item {\it Davis domains}: $2^m D_p$,  integer $m \ge 1$, for $p_{0} > p \ge 2$;
\item {\it Chebyshev-Cohn domains}: $2^m D_p$,  integer $m \ge 1$,  for $ p \ge p_{0}$;
\end{itemize}

 \begin{remark}
If the parameter $p$  is not specified $(1 < p < \infty)$, we, as in the case of Minkowski's balls, call
  Minkowski domains, Davis domains and Chebyshev-Cohn domains as Minkowski's domains.
 \end{remark}
 
\begin{proposition}
\label{p7}
Let $m$ be  integer, $m \ge 1$.
If $\Lambda$ is the critical lattice of the 
ball $D_p$ than the sublattice $\Lambda_{2^m}$  of index $2^m$ is the critical lattice of the domain $2^{m-1} D_p$.
\end{proposition}

\section{ Direct systems and Direct limits}
\subsection{ Direct systems}
Direct systems and direct limits (as well as inverse systems and inverse limits) were defined by Pontryagin \cite{Pontryagin}.
The direct system of Minkowski balls and domains has the form (\ref{dsm}), where the multiplication by $2$  is the continuous mapping
\begin{equation}
\label{dsm}
 \begin{CD}
 D_p @>2>> 2 D_p @>2>> 2^2 D_p @>2>> \cdots @>2>> 2^m D_p @>2>> \cdots 
 \end{CD}
\end{equation}

The direct system of  critical lattices  has the form (\ref{dsml}), where the multiplication by $2$  is the homomorphism of abelian groups
\begin{equation}
\label{dsml}
 \begin{CD}
 \Lambda_{p} @>2>> 2 \Lambda_{p} @>2>> 2^2 \Lambda_{p} @>2>> \cdots @>2>> 2^m \Lambda_{p} @>2>> \cdots 
 \end{CD}
\end{equation}

In our considerations we have direct systems of Minkowski balls, Minkowski domains and direct systems of critical lattices with respective maps and homomorphisms. 
Let ${\mathbb Q}_2$ and ${\mathbb Z}_2$ be respectively the field of $2$-adic numbers and its ring of integers.
Denote the corresponding direct limits by $D^{dirlim}_p$ and by $\Lambda_{p}^{dirlim}$.


\subsection{ Direct limits} 

 \begin{proposition}
  $D^{dirlim}_p = \varinjlim  2^m D_p  \in ({\mathbb Q}_2 / {\mathbb Z}_2)  D_p =
   (\bigcup_m \frac{1}{2^m} {\mathbb Z}_2/{\mathbb Z}_2) D_p.$
 \end{proposition}
 
   \begin{proposition}
  $\Lambda_{p}^{dirlim} = \varinjlim  2^m \Lambda_{p}  \in ({\mathbb Q}_2 / {\mathbb Z}_2)  \Lambda_{p} = 
  (\bigcup_m \frac{1}{2^m} {\mathbb Z}_2/{\mathbb Z}_2) \Lambda_{p}.$
 \end{proposition}
 
\begin{remark}
  In Propositions of the section $({\mathbb Q}_2 / {\mathbb Z}_2)$ is the constant  $2$-divisible group of height one.
 \end{remark}
 
\section{ Inverse systems and Inverse limits} 
\begin{remark}
   For natural $m \ge 1$ we have 
   \[
   1\cdot (2^m D_p) = \frac{1}{2}(2 \cdot 2^m D_p) =  \frac{1}{2} (2^{m+ 1} D_p) = 2^m D_p,
   \]
 \end{remark}
 and the same for critical lattices.
 \subsection{ Inverse systems}
 The inverse system of Minkowski balls and domains has the form (\ref{ism}), where the multiplication by $\frac{1}{2}$  is the continuous mapping
\begin{equation}
\label{ism}
 \begin{CD}
 D_p @<1/2<< 2 D_p @<1/2<< 2^2 D_p @<1/2<< \cdots @<1/2<< 2^m D_p @<1/2<< \cdots 
 \end{CD}
\end{equation}

The  inverse system of  critical lattices  has the form (\ref{isml}), where the multiplication by $2$  is the homomorphism of abelian groups
\begin{equation}
\label{isml}
 \begin{CD}
\Lambda_{p} @<1/2<<2 \Lambda_{p} @<1/2<< 2^2\Lambda_{p} @<1/2<< \cdots @<1/2<<2^m \Lambda_{p}@<1/2<<\cdots 
 \end{CD}
\end{equation}

 \subsection{ Inverse limits}
 
  \begin{proposition}
   The inverse (projective) limit $ D^{invlim}_p = \varprojlim 2^m D_p $ of the projective system (\ref{ism})
    is the free ${\mathbb Z}_2$-module of rank one.
   \end{proposition}
   
    \begin{proposition}
     The inverse (projective) limit 
  $\Lambda_{p}^{invlim} = \varprojlim  2^m \Lambda_{p}$
  of the projective system (\ref{isml})
   is the free ${\mathbb Z}_2$-module of rank two.
    \end{proposition}

 \end{document}